\documentclass[12pt]{article}
\usepackage{amssymb}
\usepackage{latexsym}
\usepackage{graphicx}
\usepackage{amsmath}
\usepackage{subfigure}

\newtheorem{theorem}{Theorem}

\newtheorem{cors}[theorem]{Corollary}

\newtheorem{lem}[theorem]{Lemma}

\newtheorem{observation}[theorem]{Observation}

\newtheorem{props}[theorem]{Proposition}

\newtheorem{thm}[theorem]{Theorem}

\newenvironment{pf}[1][Proof]{\noindent\textbf{#1.} }{\ \rule{0.5em}{0.5em}}
\def\lab(#1)#2{\put(#1){\makebox(0,0)[c]{#2}}}

\begin{document}

\title{A geometric characterisation of the quadratic min-power centre}

\author{M. Brazil, C.J. Ras, D.A. Thomas}
\date{}
\maketitle

\begin{abstract}
For a given set of nodes in the plane the min-power centre is a point such that the cost of the star centred at this point and spanning all nodes is minimised. The cost of the star is defined as the sum of the costs of its nodes, where the cost of a node is an increasing function of the length of its longest incident edge. The min-power centre problem provides a model for optimally locating a cluster-head amongst a set of radio transmitters, however, the problem can also be formulated within a bicriteria location model involving the $1$-centre and a generalized Fermat-Weber point, making it suitable for a variety of facility location problems. We use farthest point Voronoi diagrams and Delaunay triangulations to provide a complete geometric description of the min-power centre of a finite set of nodes in the Euclidean plane when cost is a quadratic function. This leads to a new linear-time algorithm for its construction when the convex hull of the nodes is given. We also provide an upper bound for the performance of the centroid as an approximation to the quadratic min-power centre. Finally, we briefly describe the relationship between solutions under quadratic cost and solutions under more general cost functions.\\

\noindent\textbf{Keywords}: \small{networks, power efficient range assignment, wireless ad hoc networks, generalised Fermat-Weber problem, farthest point Voronoi diagrams}
\end{abstract}

\section{Introduction}
One of the most important problems in the optimal design of wireless ad hoc radio networks is that of power minimisation. This is true during the physical design phase and when designing efficient routing protocols \cite{mon,yua,zhu}. The most appropriate fundamental model in both cases is the \textit{power efficient range assignment problem}, where communication ranges are assigned to all transmitters such that total {power} is minimised. It is assumed that the \textit{power} of each transmitter is proportional to its assigned communication range raised to an exponent $\alpha>1$ (see \cite{alt}). The process of assigning communication ranges therefore determines the total power output as well as the available communication topology of the resultant network.

\begin{figure}[htb]
    \begin{center}
        \subfigure[]{\includegraphics[scale=0.5]{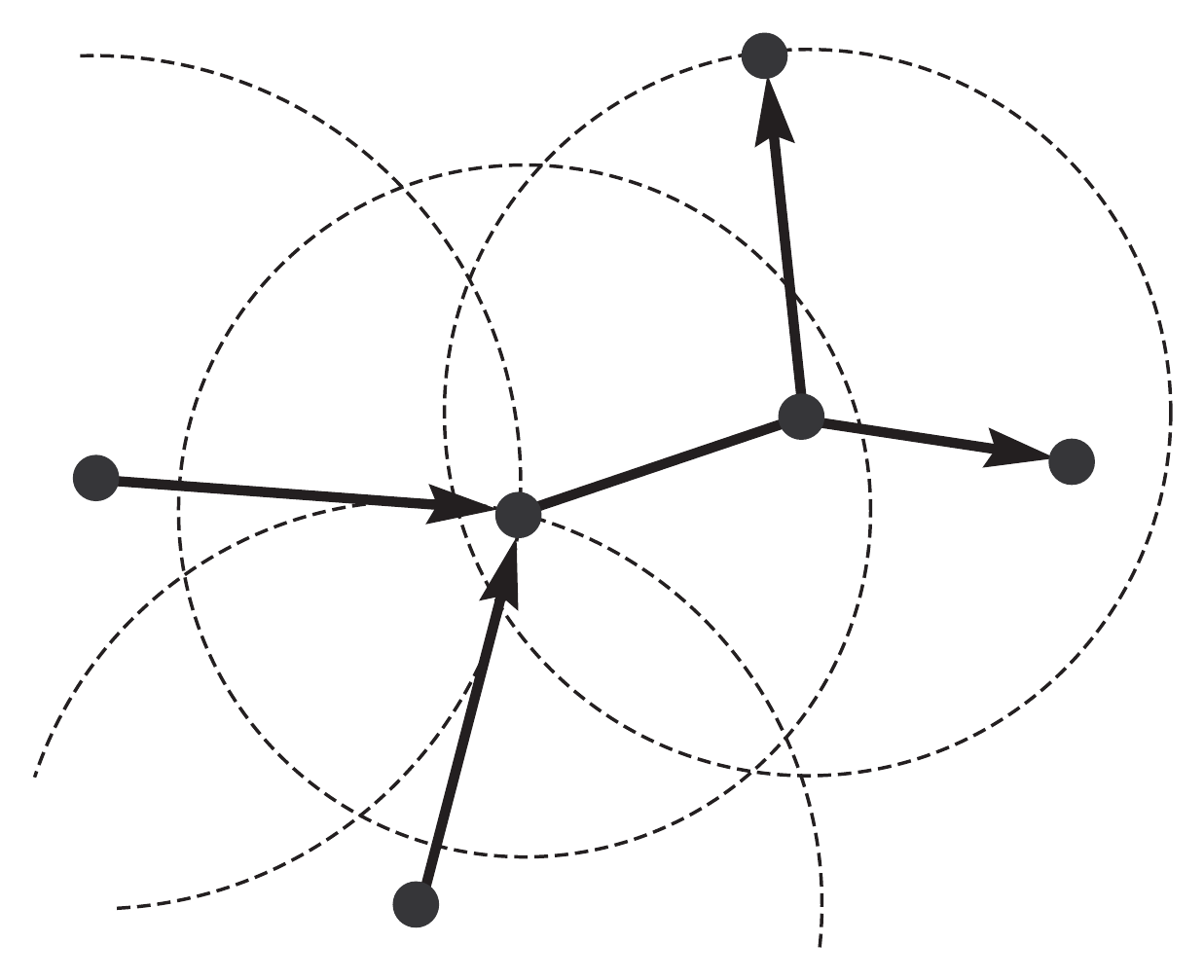}}
        \subfigure[]{\includegraphics[scale=0.5]{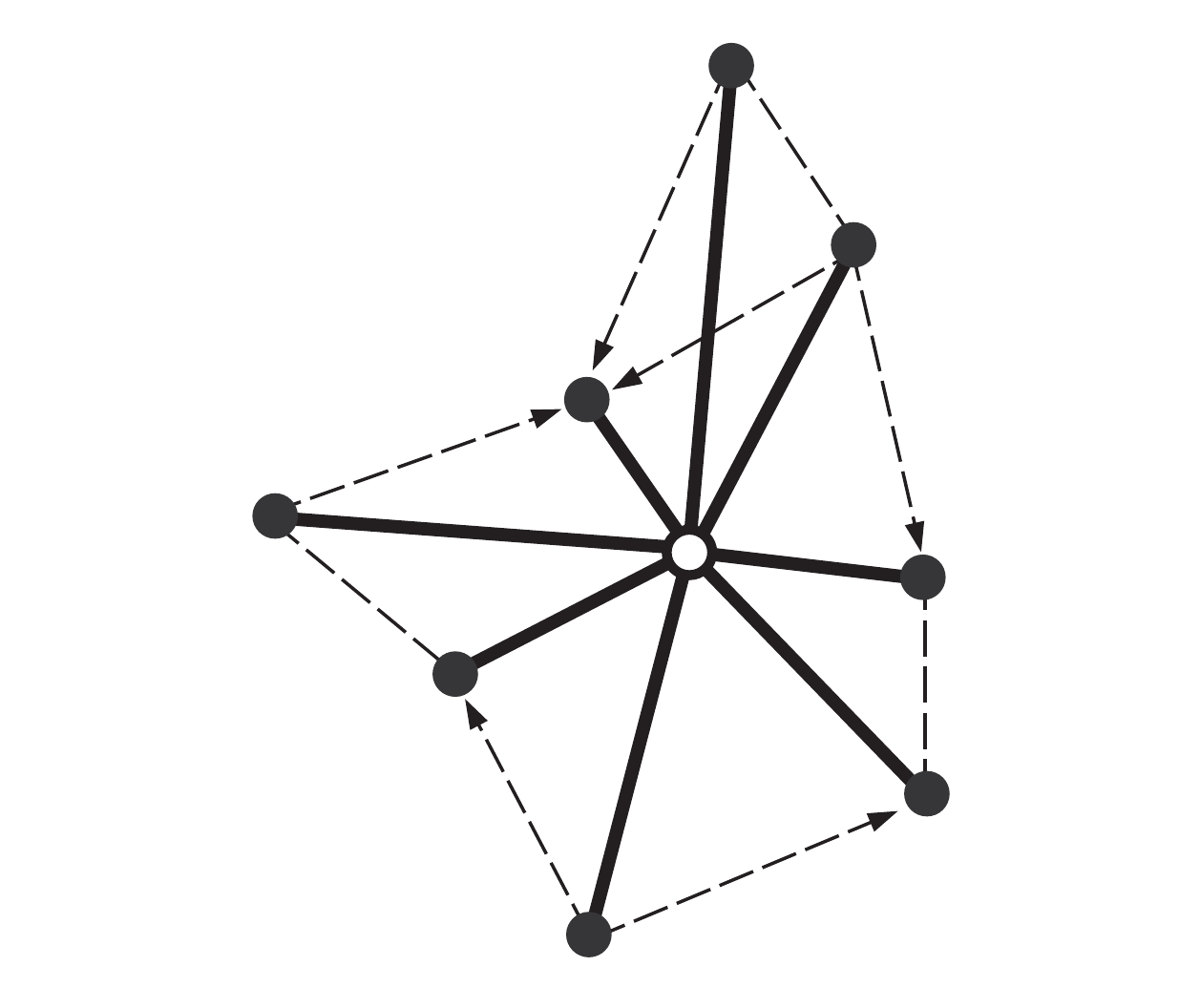}}
    \end{center}
    \caption{Examples of range assignment-induced graphs\label{figDisk}}
\end{figure}

The range assignment problem is a type of \textit{disk covering problem}, where the centres of the disks are given nodes, the radii ($r_i$) of the disks are transmission ranges, and the directed graph induced by the disks (see Figure \ref{figDisk}, where bidirected edges are depicted without arrows) must satisfy a given connectivity constraint whilst minimising $\sum r_i^\alpha$. The exponent $\alpha$ is called the \textit{path loss exponent} and most commonly takes a value between $2$ and $4$, with $\alpha=2$ corresponding to transmission in free-space. For any $\alpha>1$ the range assignment problem is NP-hard, even in the case when only $1$-connectivity is required of the resultant network \cite{fuc}.

During the design or maintenance of ad hoc radio networks it is often pertinent to introduce relays or cluster-heads for the processing of aggregated data and for the improved routing efficiency that takes place in such hierarchical structures (see \cite{dha,pau,shi}). Solving the range assignment problem whilst allowing for the introduction of a bounded number of additional nodes anywhere in the plane constitutes a very general and highly applicable geometric network problem, which has only been solved in certain restricted settings (see for instance \cite{bra1,bra2}). Since the optimal locations of the cluster-heads must be found, as well as the optimal assignment of ranges on the complete set of nodes, this so called \textit{geometric range assignment problem} is at least as difficult as the classical range assignment problem.

This paper considers the problem of optimally locating a \textit{single} cluster-head amongst a given set of transmitters, where each transmitter can send and receive data directly from the cluster-head. Not only is this an interesting and applicable model in it own right, but it is also a necessary first step in understanding the local structure of optimal networks with multiple cluster-heads. The graph induced by the assignment of ranges, in this case, contains an undirected star with the cluster head as its centre and the complete set of transmitters as its leaves; see for instance Figure \ref{figDisk}(b), where edges between transmitters are depicted as broken lines since they do not contribute to the cost (power) of the star.

In more formal terms we denote the given finite set of transmitters by $X\subset \mathbb{R}^2$ and the cluster-head by $s\in \mathbb{R}^2$. The \textit{power} of any $x\in X$ is $P_x=\|s-x\|^\alpha$ and the \textit{power} of $s$ is $P_s=\max\{\|s-x\|^\alpha:x\in X\}$, where $\|\cdot\|$ is the Euclidean norm. The \textit{total power} of the system is denoted by $P(s)=P_s+\displaystyle\sum_{x\in X}P_x$, and a \textit{min-power centre} of $X$ is a point $s^*$ which minimises $P(s)$. Minimising only $P_s$ is clearly equivalent to the $1$-\textit{centre problem}, i.e., the problem of finding the centre of a minimum spanning circle for $X$. Minimising only $\sum P_x$ is a \textit{generalised Fermat-Weber} problem \cite{bri}, which becomes the classical Fermat-Weber problem when $\alpha=1$ and the problem of constructing the \textit{centroid} when $\alpha=2$.

A related concept in facility location is the computation of the \textit{center-median} (or \textit{cent-dian}) of a finite set of points \cite{col}, which also has applications in wireless ad hoc networks for finding a so called \textit{core} node \cite{dvi}. In this problem the task is to find all Pareto-optimal solutions to the vector function $\Phi(s)=\left(P_s,\sum P_x\right)$, which is equivalent to finding the optimal value of $\lambda P_s+(1-\lambda)\sum P_x$ for every $\lambda\in [0,1]$ (see \cite{dui,fer}). Clearly the min-power centre will be optimal for this problem when $\lambda=1/2$. The vector minimisation problem has also been considered for $\alpha=1$ in the rectilinear plane \cite{mcg} and for $\alpha=2$ in the Euclidean plane \cite{ohs}. These types of bicriteria models have been described as seeking a balance between the antagonistic objectives of \textit{efficiency} (i.e., the minisum component) and \textit{equity} (i.e., the minimax component).

There are numerical methods, for instance the sub-gradient method, that optimally locate $s^*$ to within any finite precision. However, structural results for the min-power centre problem, of the type described in this paper, are necessary for optimally constructing more complex geometric range assignment networks (which is the overarching goal of our research). This fact is particularly manifest in the design of algorithmic \textit{pruning} modules, where one develops strategies based on properties of locally optimal structures for eliminating suboptimal network topologies from the exponential set of possible topologies. The ultimate benefits of good pruning modules has been demonstrated a number of times for problems similar to the geometric range assignment problem \cite{bra1,bra2,war}.

In this paper we mostly focus on the quadratic case, $\alpha=2$. In terms of cluster-head placement this means that radio transmission takes place in free space, that is, in an ideal medium with zero resistance. Path loss exponents close to $2$ frequently occur in real-world wireless radio network scenarios. This is true for transmission in mediums of low resistance, and in mediums of higher resistance when there is a degree of \textit{beam forming} (constructive interference) \cite{kar}. The quadratic case also applies to certain classical facility location problems, including the location of emergency facilities such as hospitals and fire stations \cite{fer,ohs,pue}. Furthermore, as demonstrated in the final section, it is anticipated that theoretical developments in the $\alpha=2$ case will lead to solutions and approximations for other $\alpha>1$.

In Section \ref{SecProps} we provide definitions and set up a Karush-Kuhn-Tucker formulation of the min-power centre problem and its dual for any $\alpha>1$. When $\alpha=2$ the geometric construction of the min-power centre becomes tractable, allowing us to provide a characterisation of the solution in terms of the \textit{farthest point Voronoi diagram} on $X$, and its dual, the \textit{farthest point Delaunay triangulation}. This characterisation is described in Section \ref{SecA2}, where we also present a new linear-time algorithm for the construction of the optimal quadratic min-power centre. For $\alpha=2$, Section \ref{secSignif} explores the significance of the centroid and the $1$-centre of a set of points for the min-power centre problem, and provides a characterisation of point sets for which the min-power centre and the $1$-centre coincide. Section \ref{secSignif} also provides a bound for the performance of the centroid as an approximation to the min-power centre. The final section briefly explores the general $\alpha>1$ case.

\section{Definitions and analytical properties}\label{SecProps}
Let $X=\{x_i:i\in J\}$ be a given finite set of points in the Euclidean plane with index-set $J$, and let $\alpha>1$ be a given real number. For any point $s\in \mathbb{R}^2$ let $G=G(s,X)$ be the undirected star with centre $s$ and leaf-set $X$. The \textit{power} of $G$ with respect to $s$ is denoted by
$$P(s)=P_\alpha(s,X)=\displaystyle\sum_{i\in J} \|s-{x}_i\|^\alpha+\max_{i\in J}\{\|s-{x}_i\|^\alpha\},$$ where $\|\cdot\|$ is the Euclidean norm. The first terms of of $P(s)$ can be thought of as representing the total power required by the existing transmitters for communicating with the cluster-head, while the second term can be thought of as the power required by the cluster-head for communicating with the transmitters.

\bigskip
\noindent\textbf{Definition:} The \textit{Euclidean min-power centre problem} on $X$ is the problem of finding a point $s_\alpha^*\in\underset{{s}}{\mathrm{arg\,min}}\,\, P(s)$. We refer to ${s}^*=s_\alpha^*$ as a \textit{min-power centre} of $X$.

\begin{lem}For any $X$ the function $P(s)$ is strictly convex, and therefore $s^*$ is unique.
\end{lem}
\begin{pf}
This follows since $\|s-x_i\|^\alpha$ is strictly convex when $\alpha>1$.
\end{pf}

Although $P(s)$ is not smooth, it can be expressed as the maximum of a set of smooth functions: for any $j\in J$ let $F_j(s)=F_j^\alpha(s,X)=\|s-x_j\|^{\,\alpha}+\displaystyle\sum_{i\in J} \|s-x_i\|^{\,\alpha}$. Then $P(s)=\displaystyle\max_{j\in J}\{F_j(s)\}$. Since $\alpha>1$, the gradient

\begin{displaymath}
   \nabla (\|s-x_i\|^\alpha)= \left\{
     \begin{array}{lr}
       \alpha (s-x_i)\|s-x_i\|^{\alpha-2} &: s\neq x_i,\\
       0  &: s=x_i,
     \end{array}
   \right.
\end{displaymath}
is continuous. Therefore $F_j$ is continuously differentiable and convex, since it is the sum of such functions.

Next we describe a useful non-linear programming formulation of the problem of finding $s^*$. Consider the following non-linear optimization program.\\

\noindent\textbf{Problem P1}.
\begin{align}
&\underset{{(s,v)\in\mathbb{R}^2\times\mathbb{R}}}{\mathrm{minimise}}\,\,v\nonumber\\
&\mathrm{subject\ to\ }\nonumber\\
&F_j(s)-v\leq 0,\ \forall j\in J.\nonumber
\end{align}

Clearly Problem P1 is smooth and convex. It is easy to see that $(s,v)$ is optimal for this problem if and only if $s=s^*$ and $v=P(s^*)$.

\begin{lem}Problem P1 satisfies the Mangasarian-Fromovitz constraint qualification for any point $(s,v)\in \mathbb{R}^2\times\mathbb{R}$.
\end{lem}
\begin{pf}
Note that $$\langle\nabla (F_j(s)-v),((0,0),1)\rangle=\langle(\nabla F_j(s),-1),((0,0),1)\rangle<0$$ for all $j$ such that $P(s)=F_j(s)$.
\end{pf}

By the previous lemma and the resulting KKT conditions it follows that the point $(s^*,v^*)\in \mathbb{R}^2\times\mathbb{R}$ is optimal for Problem P1 if and only if there exist multipliers $\lambda_j^*\geq 0,\ j\in J$ such that $\nabla v+\displaystyle\sum_{j\in J}\lambda_j^*(\nabla F_j(s^*),-1)=0$ and $\lambda_j^*(F_j(s^*)-v^*)=0$ for all $j\in J$. Observe that

\begin{alignat*}{5}
&\lambda_j^*(F_j(s^*)-v^*)=0& &\Leftrightarrow \lambda_j^*(F_j(s^*)-P(s^*))=0&\\
&{}& &\Leftrightarrow\lambda^*_j(\|s^*-x_j\|^\alpha-\max_{i\in J}\{\|s^*-x_i\|^\alpha\})=0&\\
&{}& &\Leftrightarrow\lambda^*_j(\|s^*-x_j\|-\max_{i\in J}\{\|s^*-x_i\|\})=0&
\end{alignat*}

Therefore we can state equivalent KKT conditions as follows: $s^*$ is the min-power centre of $X$ if and only if there exists $\lambda^*_j\geq 0$ such that

\begin{align}
      &\displaystyle\sum_{j\in J}\lambda_j^* \nabla F_j(s^*) =0,\label{eq:KKT1a}\\
      &\displaystyle\sum_{j\in J}\lambda_j^* =1,\label{eq:KKT2a}\\
      &\lambda^*_j(\|s^*-x_j\|-\max_{i\in J}\{\|s^*-x_i\|\})=0,\quad \forall j\in J.\label{eq:KKT3a}
\end{align}

Let $n=|J|$ and let $\lambda=(\lambda_1,...,\lambda_n)$. The dual of Problem P1 can be written as\\

\noindent\textbf{Problem  D1.}
\begin{align}
&\underset{{\lambda\in\mathbb{R}^n}}{\mathrm{maximise}}\,\,\underset{{s\in\mathbb{R}^2}}{\min}
\displaystyle\sum_{j\in J}(1+\lambda_j)\|s-x_j\|^\alpha\nonumber\\
&\mathrm{subject\ to\ }\nonumber\\
&\lambda \geq 0,\nonumber\\
&\displaystyle\sum_{j\in J}\lambda_j=1.\nonumber
\end{align}

Finally, we can convert the dual problem to a more familiar form as follows. For every $j\in J$ let $\tau_j=\dfrac{\lambda_j+1}{n+1}$. Therefore $\dfrac{1}{n+1}\leq \tau_j$ and $\displaystyle\sum_{j\in J} \tau_j=1$. Let $s(\tau)=\underset{{s\in\mathbb{R}^2}}{\mathrm{arg}\,\min}\displaystyle\sum_{j\in J}\tau_j\|s-x_j\|^\alpha$ (observe that $s(\tau)$ is unique). Clearly then $\displaystyle\sum_{j\in J}\tau_j\nabla(\|s(\tau)-x_j\|^\alpha)=0$, which is equivalent to $\displaystyle\sum_{j\in J}\tau_j(s(\tau)-x_j)\|s(\tau)-x_j\|^{\alpha-2}=0$. Therefore an equivalent formulation of Problem D1 is:\\

\textbf{Problem D1$'$.}
\begin{align}
&\underset{{\tau\in\mathbb{R}^n}}{\mathrm{maximise}}\,\,(n+1)\displaystyle\sum_{j\in J}\tau_j\|s-x_j\|^\alpha\label{eq:D1}\\
&\mathrm{subject\ to\ }\nonumber\\
&\dfrac{1}{n+1}\leq \tau_j,\label{eq:D2}\\
&\displaystyle\sum_{j\in J}\tau_j=1,\label{eq:D3}\\
&\displaystyle\sum_{j\in J}\tau_j(s-x_j)\|s-x_j\|^{\alpha-2}=0.\label{eq:D4}
\end{align}

If $\tau^*$ is optimal for Problem D1$'$ then $s^*=s(\tau^*)$ is optimal for Problem P1 and $\lambda^*=(n+1)\tau^*-1$ gives the multipliers satisfying Equations (\ref{eq:KKT1a})-(\ref{eq:KKT3a}). This final formulation is recognisable in that it is almost identical to the dual of the problem $$\underset{s\in\mathbb{R}^2}{\mathrm{minimise}}\displaystyle\max_{j\in J}||s-x_j||^\alpha,$$ which is the well-known \textit{minimum spanning circle problem} (see \cite{elz} for a derivation). The only two differences in these duals are Inequality (\ref{eq:D2}), which will be replaced by $\tau_j\geq 0$ for the minimum spanning circle problem, and the extra coefficient $n+1$ of the objective function. In general, the restriction on domain specified by Inequality \ref{eq:D2} leads to a non-trivial variation on the minimum spanning circle problem, however, as shown in the next section, for $\alpha=2$ there exists an interesting and useful geometric relationship between the minimum spanning circle problem and the min-power centre problem.

\section{The geometry of $\alpha=2$}\label{SecA2}
In order to study the geometry of the case where $\alpha=2$ we first rewrite the KKT conditions in terms of the $2$-\textit{centroids} of $X$.\\

\noindent\textbf{Definition:} Recall that the \textit{centroid} (or centre of mass) of $X$, which we denote as $M$, is defined by $M=\dfrac{1}{n}\displaystyle\sum_{j\in J}x_j$. We define the set of $2$-\textit{centroids} of $X$ to be the set $\mathcal{M}=\{M_j:j\in J\}$ where each $M_j=\dfrac{1}{n+1}\left(x_j+\displaystyle\sum_{i\in J}x_i\right)$.

\begin{figure}[htb]
  \begin{center}
    \includegraphics[scale=0.45]{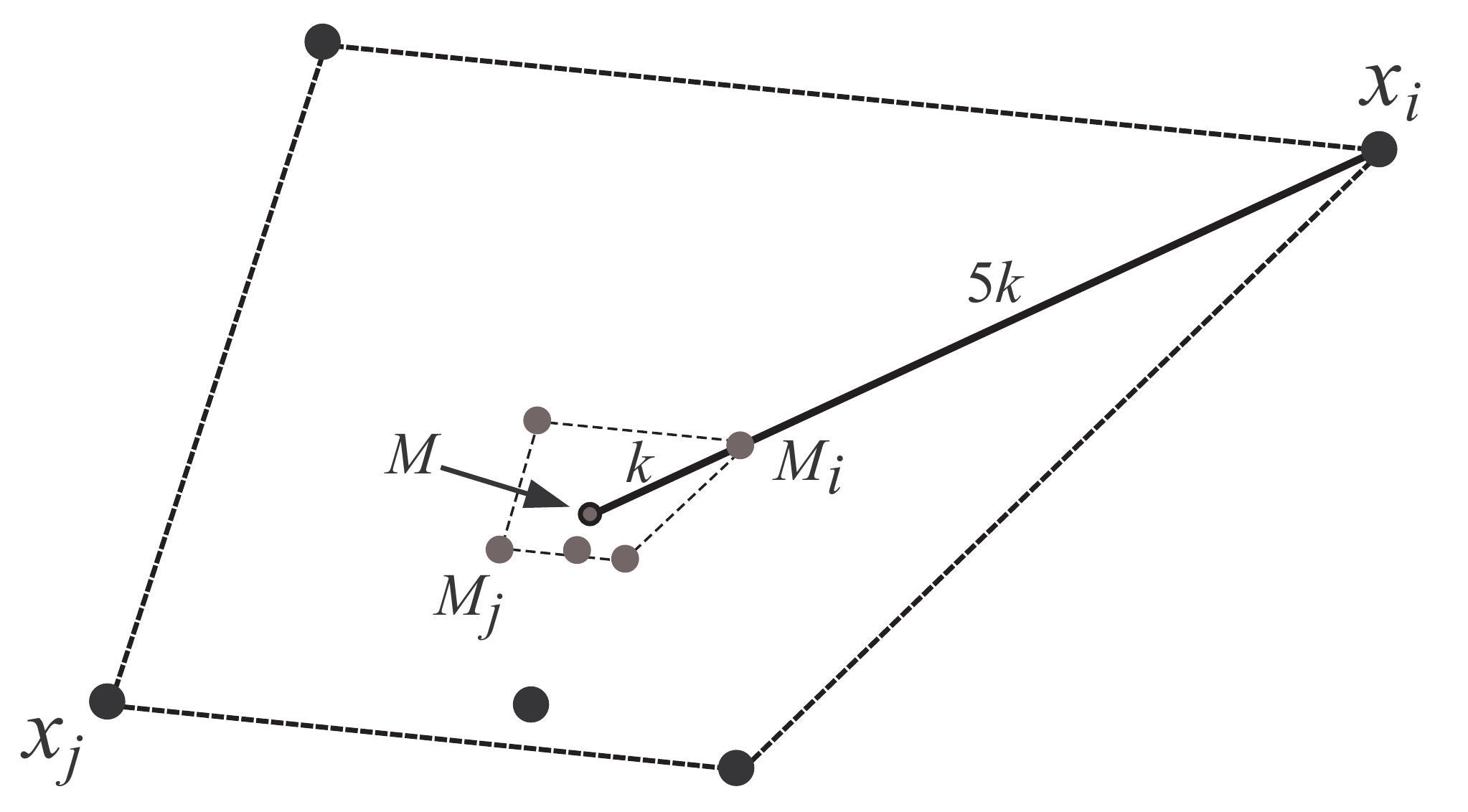}\\
  \end{center}
  \caption{An example of a construction of $2$-centroids $\mathcal{M}=\{M_j\}$ for $n=5$. For every $j$ we have $||M-M_j||:||M_j-x_j||=1:n$}
  \label{figEx1}
\end{figure}

Observe that (if $X$ and $\mathcal{M}$ are represented as vectors) $\mathcal{M}=AX$ where $A=\frac{1}{n+1}(I+1_n)$, $I$ is the $n\times n$ identity matrix, and $1_n=(b_{ij})$ is the $n\times n$ matrix with $b_{ij}=1$ for all $i,j$. This means that $\mathcal{M}$ is the image of an affine transformation on $X$. An consequence of this observation is the following lemma. Let $\mathrm{conv}(\cdot)$ denote the the convex hull of a set of points.

\begin{lem}\label{divrat}The region $\mathrm{conv}(\mathcal{M})$ is similar to $\mathrm{conv}(X)$, and $M$ is the centroid of $\mathcal{M}$. For any $j$ the point $M_j$ divides the line segment $Mx_j$ into a ratio of $1:n$ (see Figure \ref{figEx1}, where the $M_j$ are grey-filled circles).
\end{lem}

For any $X'\subseteq \mathrm{conv}(X)$ we denote (with a slight abuse of notation) the image of $X'$ under the above affine transformation by $AX'$, which is a subset of $\mathrm{conv}(\mathcal{M})$. We also similarly employ $A^{-1}$ for the inverse transformation.

In the quadratic case (where $\alpha=2$) we now write
\begin{alignat*}{5}
&\nabla F_j(s)& &=2(n+1)s-2x_j-2\displaystyle\sum_{i\in J}x_i&\\
&{}& &=2(n+1)(s - M_j).&
\end{alignat*}

This equation implies that for any $j$ we have $M_j=\underset{{s}}{\mathrm{arg\,min}}\,\, \{\|x_j-s\|^2+\sum_{i\in J} \|x_i-s\|^2\}$. In the case $\alpha=2$ a simplified set of KKT conditions for Problem P1 in total is
\begin{align}
      &s^*=\displaystyle\sum_{j\in J}\lambda_j^* M_j,\label{eq:KKT1}\\
      &\displaystyle\sum_{j\in J}\lambda_j^*=1,\label{eq:KKT2}\\
      &\lambda^*_j(\|s^*-x_j\|-\max_{i\in J}\{\|s^*-x_i\|\})=0,\quad \forall j\in J.\label{eq:KKT3}
\end{align}

By Condition \ref{eq:KKT1} we have:
\begin{lem}$s^*\in\mathrm{conv}(\mathcal{M})$.
\end{lem}

Conditions (\ref{eq:KKT1})--(\ref{eq:KKT3}) can be interpreted by means of a dual pair of planar geometric graphs. Firstly let $\mathcal{V}=\mathcal{V}(X)$ be the \textit{farthest point Voronoi diagram} on $X$. This is a partition of $\mathbb{R}^2$ into maximal regions $\{V(x_j):j\in J\}$ such that $s\in V(x_j)$ if and only $x_j$ is a farthest node from $s$. It is well-known that $V(x)$ is non-empty if and only if $x$ is an \textit{extreme point} (i.e., a ``corner") of $\mathrm{conv}(X)$. As is common practice, we will (without confusion) treat $\mathcal{V}$ both as a graph specified by the boundary edges, or as a partition of the plane. It is well known that $\mathcal{V}$ is a tree.

\begin{figure}[htb]
  \begin{center}
    \includegraphics*[scale=0.6, trim= 50 80 50 70]{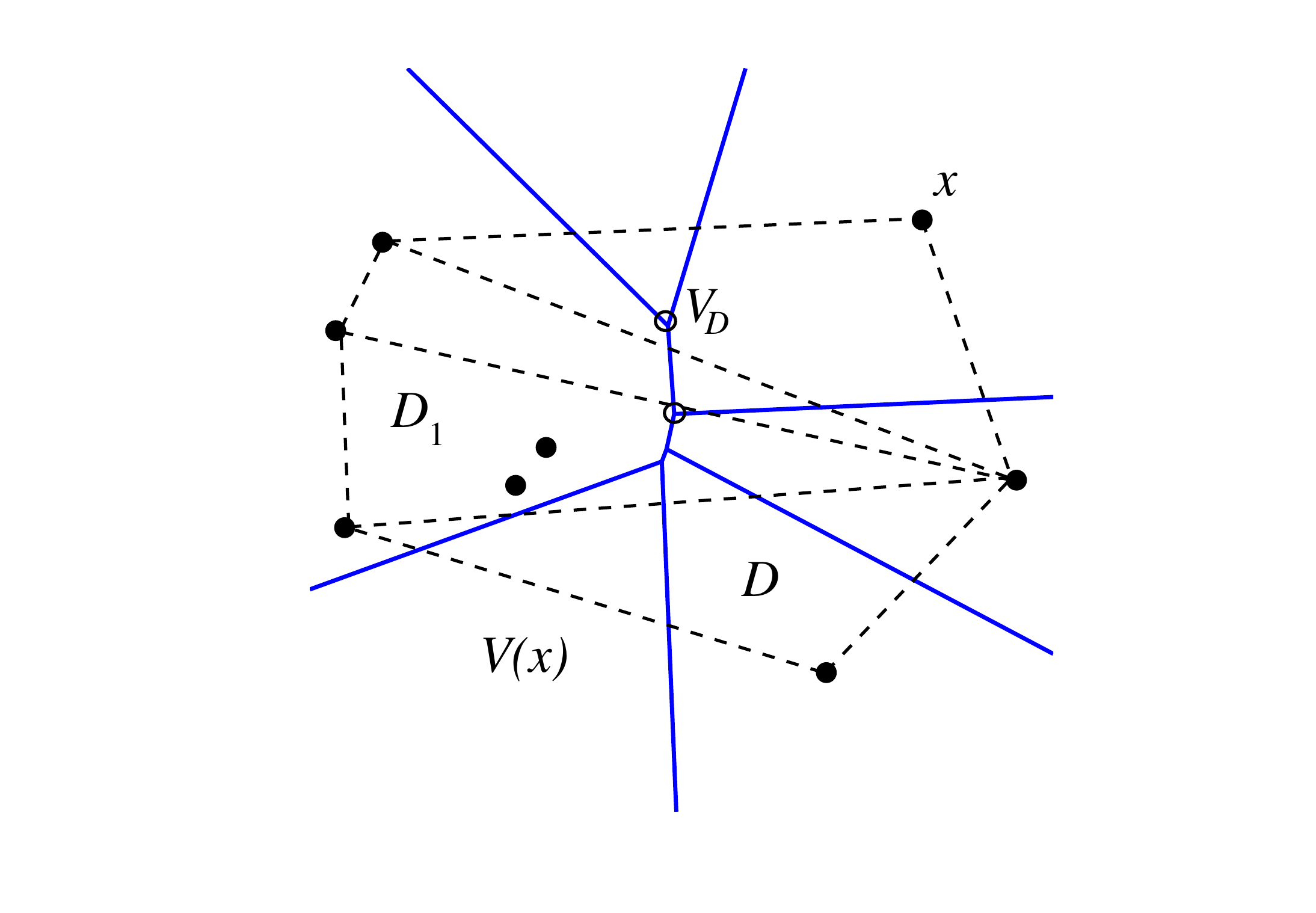}\\
  \end{center}
  \caption{A farthest point Voronoi diagram (solid lines) and its dual (broken lines) on a set of eight nodes (black-filled circles)}
  \label{figEx2}
\end{figure}

Let $\mathcal{D}=\mathcal{D}(X)$ be a \textit{farthest point Delaunay triangulation} of $X$, which is defined as a triangulation of $X$ such that the circumcircle of any triangle contains all points of $X$. The farthest point Delaunay triangulation is the dual graph of $\mathcal{V}$ \cite{epp}. If $\mathcal{V}$ contains a vertex of degree more than $3$ then $\mathcal{D}$ is not unique; in this case we let $\mathcal{D}$ be \textit{any} farthest point Delaunay triangulation on $X$. Figure \ref{figEx2} illustrates an example of both these graphs.

A $b$-\textit{face} of $\mathcal{V}$ or $\mathcal{D}$ is a maximal region of dimension $b$, where $0\leq b\leq 2$. Therefore a $0$-face is a vertex, a $1$-face is an edge, and a $2$-face is a region of $\mathcal{V}$ or $\mathcal{D}$. The dual relationship between $\mathcal{D}$ and $\mathcal{V}$ defines a natural mapping from each $b$-face of $\mathcal{D}$ to a $(2-b)$-face of $\mathcal{V}$. In Figure \ref{figEx2} the $0$-face $x$ of $\mathcal{D}$ maps to the $2$-face $V(x)$ of $\mathcal{V}$, and the $2$-face $D$ of $\mathcal{D}$ maps to the $0$-face $V_D$ of $\mathcal{V}$. Observe that dual $1$-faces are always perpendicular to each other. If the maximum degree of $\mathcal{V}$ is at most $3$ then the mapping between faces is a bijection. In the case of vertices of degree more than $3$ we still get a bijection if we consider every such vertex as consisting of a number of degree $3$ vertices connected by zero-length edges, where the topology is uniquely determined by the choice of $\mathcal{D}$. The unique pairs in the bijection will be referred to as \textit{dual faces}. For any $b$-face $D$ of $\mathcal{D}$ we denote by $V_D$ its dual face. Formally the bijection is defined as $D\leftrightarrow V_D$ if and only if $V_D=\bigcap_{j\in J_D}V(x_j)$ where $J_D$ is the index set of the extreme points of $D$.

Observe in Figure \ref{figEx2} that there exists a Delaunay triangle $D_1$ that strictly contains its dual face in $\mathcal{V}$, depicted by the unlabelled empty circle. As will become clear below, this dual relationship is unique and, in fact, the unlabelled empty circle is the centre of the minimum spanning circle (called the $1$-\textit{centre}) of $X$. We will show that something similar holds for the min-power centre problem.

Before proving the main theorem, we note that it is often the case when solving minimax non-linear programs that a complete solution can be arrived at by considering the various combinations of functions which are potentially active at the optimal solution. Hence we introduce the following lemma which allows us to restrict the number of combinations.

\begin{lem}There exists a set of multipliers $\{\lambda_j:j\in J\}$ satisfying KKT Conditions (\ref{eq:KKT1})-(\ref{eq:KKT3}) and containing at most $3$ non-zero elements.
\end{lem}
\begin{pf}
Suppose that $s^*=\displaystyle\sum_{j\in J}\lambda_j^* M_j$ and that $J^0=\{j:\lambda_j^*>0,j\in J\}$ contains at least $4$ elements. Therefore $s^*=\displaystyle\sum_{j\in J^0}\lambda_j^* M_j$ and $s^*\in \mathrm{conv}(\mathcal{M}^0)$, where $\mathcal{M}^0=\{M_j:j\in J^0\}$. For any $j\in J^0$, since $F_j(s^*)=P(s^*)$ it follows that $x_j$ lies on a vertex of $\mathcal{V}$ and is an extreme point of $\mathrm{conv}(X)$. Therefore, by Lemma \ref{divrat}, $M_j$ is an extreme point of $\mathrm{conv}({\mathcal{M}})$ for all $j\in J^0$. By Carath\'eodory's theorem any point in the convex hull of $\mathcal{M}^0$ can be expressed as a convex combination of at most three extreme points. Let $s^*=\displaystyle\sum_{j\in J^*}u_j M_j$ for some $J^*\subset J^0$ containing at most three elements, where $u_j>0$ and $\displaystyle\sum_{j\in J^*}u_j=1$. Let $\lambda_j=0$ for all $j\in J\backslash J^*$ and let $\lambda_j=u_j$ for $j\in J^*$. Then $\{\lambda_j:j\in J\}$ and $s^*$ satisfy the KKT conditions.
\end{pf}

From now on it is assumed that the KKT set associated with $s^*$ is $\{\lambda_j^*:j\in J\}=\Lambda\cup \{\lambda_j^*:j\in J\backslash J^*\}$, where $\lambda_j^*>0$ if and only if $j\in J^*$, and $\Lambda=\{\lambda^*_j:j\in J^*\}$ is of minimal cardinality. From the above proof it follows that there may be more than one choice for this minimal set of multipliers. For consistency we assume that the multipliers are chosen so that there exists a Delaunay face $D\in \mathcal{D}$ such that $\{x_j:j\in J^*\}$ is the set of extreme points of $D$.

Let $\mathrm{int}(\cdot)$ denote the interior of a closed region ($b$-face). For ease of notation we assume that $\mathrm{int}(\{x\})=x$ for any point $x$.

\begin{figure}[!p]
    \begin{center}
        \subfigure[$|\Lambda|=1$]{\includegraphics*[scale=0.49, trim= 100 30 100 30]{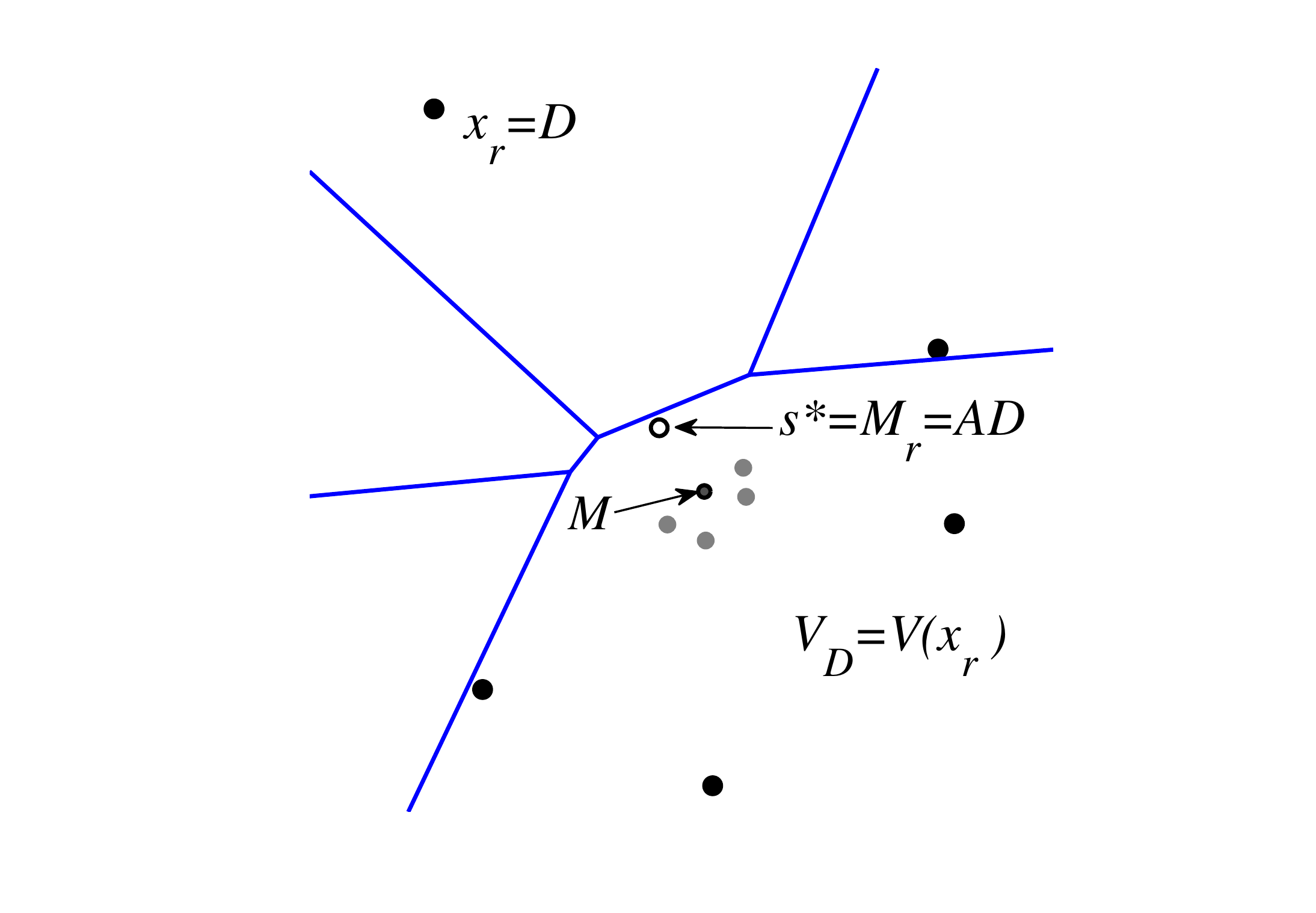}}
        \subfigure[$|\Lambda|=2$]{\includegraphics*[scale=0.49, trim= 120 20 120 20]{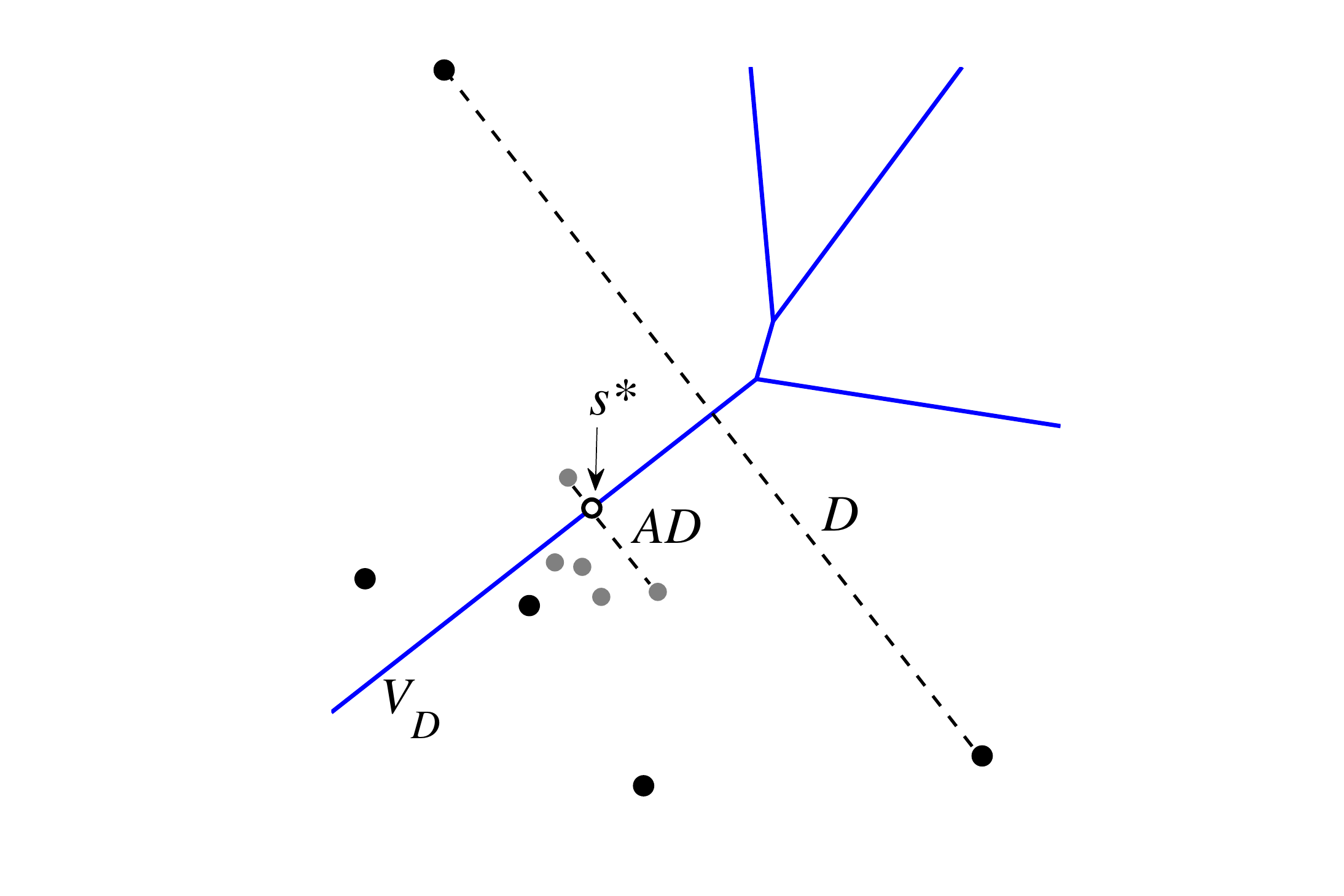}}
        \subfigure[$|\Lambda|=3$]{\includegraphics*[scale=0.49, trim= 120 20 120 20]{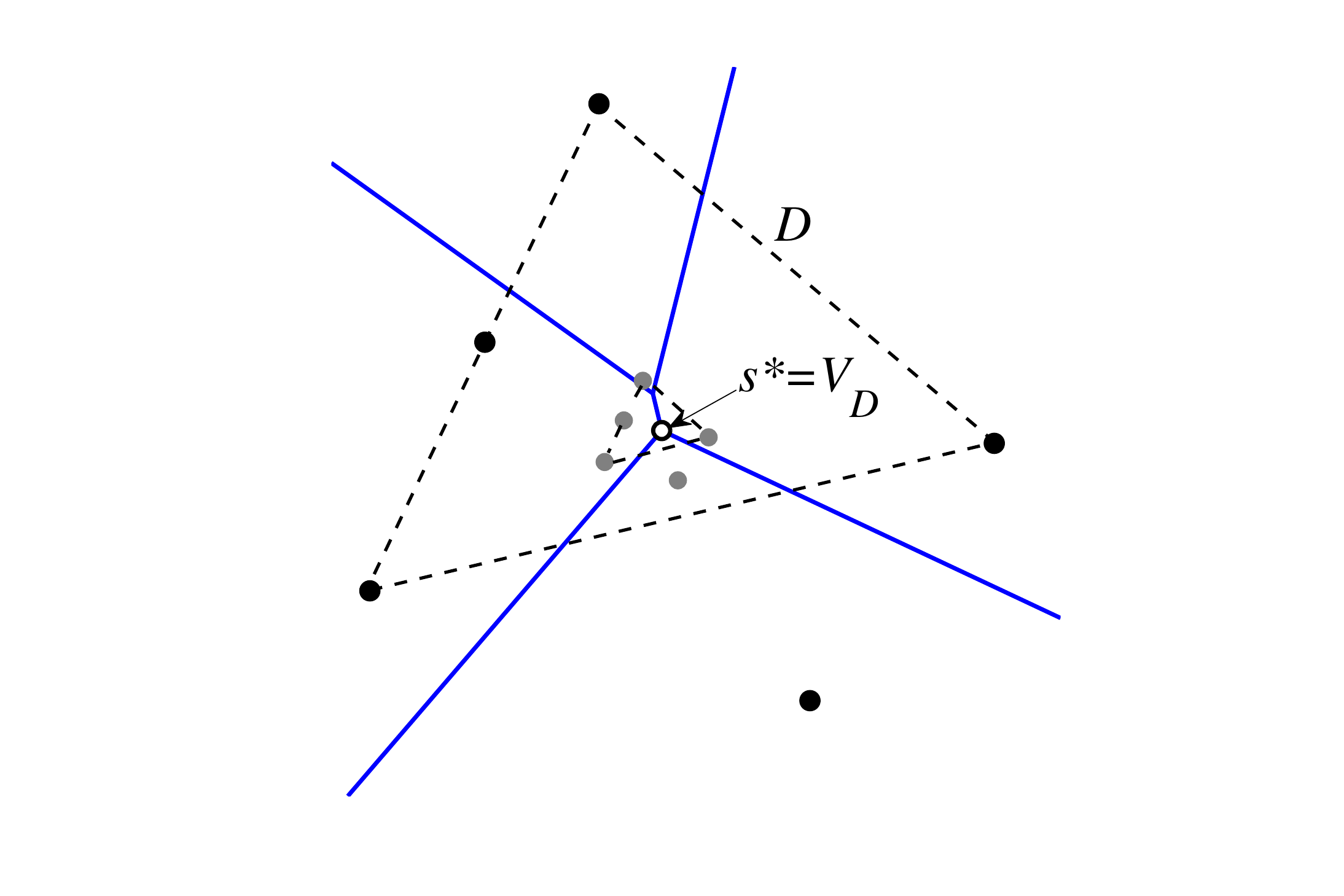}}
    \end{center}
    \caption{Illustrating the three cases of Theorem \ref{mainT} with $|J|=5$. The $x_j$ are black-filled circles, the $M_j$ grey-filled, and $s^*$ white-filled\label{figBotDef}}
\end{figure}

\begin{thm}\label{mainT}A point $s$ is the optimal min-power centre for $X$ if and only if $s\in \mathrm{int}(AD)$ and $s\in V_D$ for some face $D$ of $\mathcal{D}$.
\end{thm}
\begin{pf}
Suppose that $s=s^*$ is optimal and let $R$ be the face (node, edge, or triangle) defined by $R=\mathrm{conv}(\{M_j:j\in J^*\})$. By complementary slackness we have $\|s^*-x_j\|=\displaystyle \max_{i\in J}\{\|s^*-x_i\|\}$ for every $j\in J^*$. Therefore the disk centred at $s^*$ and containing the points $\{x_j:j\in J^*\}$ on its boundary covers all points of $X$. By the definition of farthest point Delaunay triangulations, $A^{-1}R=\mathrm{conv}(\{x_j:j\in J^*\})$ is a face of $\mathcal{D}$. Since $J^*$ has minimal cardinality and $s^*$ is a convex combination of the extreme points of $R$ we must have $s^*\in \mathrm{int}(R)$. Finally, note that $s^*\in\bigcap_{j\in J^*} V(x_j)$, which implies that $s^*$ is also an element of the dual face of $A^{-1}R$. Letting $D=A^{-1}R$ proves necessity.

To prove sufficiency we assume that $s\in \mathrm{int}(AD)$ and $s\in V_D$ for some face $D$ of $\mathcal{D}$. Let $J_D=\{j:M_j\mathrm{\ is\ an\ extreme\ point\ of\ }AD\}$ and let $\{u_j:j\in J_D\}$ be a set of positive weights such that $\sum u_j=1$ and $s=\sum u_jM_j$. Let $\{\lambda_j:j\in J\}$ be a set of multipliers such that $\lambda_i=u_i$ for all $i\in J_D$, and all other $\lambda_i=0$. We only need to show still that complementary slackness holds with respect to $s,\lambda$. Let $x_k$ be any point of $X$. Since $V_D$ is the dual of $D$ we have $s\in\bigcap_{j\in J_D} V(x_j)$. Therefore if $k\in J_D$ then $\|s-x_k\|-\max\{\|s-x_i\|:i\in J\}=0$. Therefore $\lambda_j(\|s-x_j\|-\max\{\|s-x_i\|:i\in J\})=0$ for all $j\in J$, which implies that $s$ is optimal.
\end{pf}

Figures \ref{figBotDef} illustrates the three different cases that arise based on the number of active functions at $s^*$. Something more can be said about the case $|\Lambda|=1$:

\begin{props}Let $x_r$ be the node of $X$ farthest from $M$. The following statements are equivalent:
\begin{enumerate}
    \item $|\Lambda|=1$,
    \item $s^*=M_r$,
    \item $M_j\in V(x_j)$ for some $j\in J$.
\end{enumerate}
\end{props}
\begin{pf}
Assume that $|\Lambda|=1$. Then, by Theorem \ref{mainT}, $s^*=M_i$ and $s^*\in V(x_i)$ for some $i$. Since $x_i$ is the farthest node from $s^*$, and $x_i,M_i,M$ are collinear, it follows that $M\in \mathrm{int}(V(x_i))$. Therefore $i=r$. Finally, suppose that $M_j\in V(x_j)$ for some $j\in J$. Then, also by Theorem \ref{mainT}, $M_j=s^*$. Therefore $|\Lambda|=1$.
\end{pf}

\begin{cors}If $s^*\in \mathrm{int}(V(x_j))$ for some $j\in J$, then $s^*=M_r$.
\end{cors}

The convention of using $x_r$ to represent the node of $X$ farthest from $M$ will hold throughout the rest of this paper. Using Theorem \ref{mainT} it is possible to design a fast geometric algorithm for finding $s^*$.\\

\noindent\textbf{Algorithm}\\
\noindent\textbf{Input}{: the set of coordinates $X$}\\
\noindent\textbf{Output}{: the point $s^*$}\\
\texttt{-\hspace{0.5cm}Construct $\mathcal{V}$}\\
\texttt{-\hspace{0.5cm}Construct $M_r$. Perform point location for $M_r$ in $V(x_r)$}\\
\texttt{-\hspace{0.5cm}If $M_r\in V(x_r)$ then $s^*=M_r$ and exit}\\
\texttt{-\hspace{0.5cm}For every edge $e=V(x_i)\cap V(x_j)$ do}\\
\texttt{-\hspace{1.5cm}Construct $M_iM_j$}\\
\texttt{-\hspace{1.5cm}If $M_iM_j\cap e\neq \emptyset$ then $s^*=M_iM_j\cap e$ and exit}\\
\texttt{-\hspace{0.5cm}For every vertex $v=V(x_i)\cap V(x_j)\cap V(x_k)$ do}\\
\texttt{-\hspace{1.5cm}Construct $\triangle M_iM_jM_k$}\\
\texttt{-\hspace{1.5cm}Perform point location for $v$ in $\triangle M_iM_jM_k$}\\
\texttt{-\hspace{1.5cm}If $v\in \triangle M_iM_jM_k$} then $s^*=v$ and exit\\

Constructing the farthest point Voronoi diagram can be done in $O(n\log n)$ time from scratch or in $O(n)$ time if the convex hull of $X$ is given. From any face $V_D$ of $\mathcal{V}$ the face $AD$ can be constructed in constant time. The point location for $M_r$ in $V(x_r)$ requires $O(\log n)$ time, but the point location for any $v$ in a $\triangle M_iM_jM_k$ requires only constant time. Since the structural complexity of $\mathcal{V}$ is $O(n)$, analysing every face of $\mathcal{V}$ until Theorem \ref{mainT} is confirmed also takes $O(n)$ time. The total complexity of the algorithm is therefore $O(n)$ when the convex hull of $X$ is given.

Ohsawa's algorithm \cite{ohs} for finding the efficient solutions to the vector minimisation problem $\mathrm{v}$-$\min\Phi(s)$, where $\Phi(s)=(\max\{\|s-x_j\|\},F(s))$, can also find $s^*$. There is no theoretical gain in the runtime of our algorithm over the algorithm of Ohsawa, but with modifications our algorithm should perform better in practice. We provide an intuitive reason for this as follows. As is well-known, for any set of $n$ points $X$ the expected number of points that are extreme points of $\mathrm{conv}(X)$ is $O(\log n)$. Therefore the expected complexity of the farthest point Voronoi diagram does not increase very quickly as $n$ increases. Now both $M$ and $s^*$ are in $\mathrm{conv}(\mathcal{M})$, but $\mathrm{conv}(\mathcal{M})$ becomes smaller the more we increase $n$ (note that $\mathcal{M}$ is contained in a disk of radius $\frac{1}{n+1}\|M-x_r\|$). Therefore we suggest that if we modify our algorithm so that it begins by checking edges and vertices of $\mathcal{V}$ close to $M$, that the expected number of steps before $s^*$ is found will be small for large $n$. In contrast, for a fixed polygon representing a convex hull, and a fixed point in the polygon representing a centroid, the cardinality of a set $X$ with this convex hull and centroid will not affect the practical runtime of Ohsawa's algorithm. Of course, Ohsawa's algorithm could also be modified to begin with points close to $M$, but his paper presents no theoretical justification for doing this.

\section{The significance of the centroid and the $1$-centre to the min-power centre problem}\label{secSignif}
It is clear from the results so far that the centroid $M$ plays a central role in the min-power centre problem when $\alpha=2$. In particular, it is intuitive that $M$ is always relatively close to the min-power centre, $s^*$, since both lie in $\mathrm{conv}(\mathcal{M})$. This fact leads to a theorem at the end of this section which places an upper bound on the performance ratio of approximating the min-power centre by the centroid. The fact that $M$ is a ``special" point in the min-power centre problem was also observed in \cite{ohs}, where it is shown that any Parato optimal solution to the vector minimisation problem $\mathrm{v}$-$\min\Phi(s)$ lies on a simple path connecting $M$ and the $1$-centre of $X$ (recall that $\Phi(s)=(P_s,F(s))$, where $P_s=\max\{\|s-x\|^2:x\in X\}$ and $F(s)=\displaystyle\sum_{x\in X}\|s-x\|^2$). Interestingly, this path is piecewise linear, consisting of one segment joining $M$ to a point on an edge of $\mathcal{V}$, and with the remainder of the path lying on the edges of $\mathcal{V}$. Since the min-power centre is Pareto optimal for $\mathrm{v}$-$\min\Phi(s)$ it follows that $s^*$ is also on this path. This brings us to a second question: what role does the $1$-centre play in the min-power centre problem? We will address this question next before returning to the discussion on $M$.

Denote the $1$-centre of $X$ by $C$ and let $\mathcal{M}^*=\{M_j:j\in J^*\}$, where, as before, $J^*$ is minimal. Therefore $s^*\in\mathrm{int}(\mathrm{conv}(\mathcal{M}^*))$. The next lemma is given without proof as it follows from an argument similar to the proof of Theorem \ref{mainT}.

\begin{lem}\label{onemore}A point $s$ is the $1$-centre of $X$ if and only if $s\in \mathrm{int}(D)$ and $s\in V_D$ for some face $D$ of $\mathcal{D}$.
\end{lem}

\begin{props}If $M\in\mathrm{conv}(\mathcal{M}^*)$ then $s^*=C$.
\end{props}
\begin{pf}
Clearly $|J^*|>1$, since if $|J^*|=1$ then $s^*=M_r$, but $M_r$ and $M$ are always distinct. Let $X^*=A^{-1}\mathcal{M}^*$. Since $M\in \mathrm{conv}(\mathcal{M}^*)$ we have, by Lemma \ref{divrat}, the fact that $\mathrm{conv}(\mathcal{M}^*)$ does not intersect the exterior of $\mathrm{conv}(X^*)$. Therefore, since $s^*\in\mathrm{int}(\mathrm{conv}(\mathcal{M}^*))$ it follows that $s^*\in\mathrm{int}(\mathrm{conv}(X^*))$. By letting $D=\mathrm{conv}(X^*)$ it follows from Lemma \ref{onemore} that $s^*=C$ as required.
\end{pf}

The converse of this result does not appear to hold, but we do have an interesting corollary.

\begin{cors}$M$ is the min-power centre of $X$ if and only if $M$ and $C$ coincide.
\end{cors}

Our final result on the role of $C$ in the min-power centre problem relates to point sets $X$ that are more or less evenly distributed (with respect to distance) about their centroid. The proof utilises the formulation of Problem D1$'$, which, as mentioned before, is almost exactly the dual of the minimum spanning circle problem except that the domain of the $\tau_j$ is further restricted by Inequality (\ref{eq:D2}). Note that when $\alpha=2$ we can write Equation (\ref{eq:D4}) as $s=\sum\tau_jx_j$ since $\sum\tau_j=1$. Therefore, for $\alpha=2$ we have:\\

\textbf{Problem D2 (dual of the minimum spanning circle problem).}
\begin{align}
&\underset{{\tau\in\mathbb{R}^n}}{\mathrm{maximise}}\,\,\displaystyle\sum_{j\in J}\tau_j\|s(\tau)-x_j\|^2\nonumber\\
&\mathrm{subject\ to\ }\nonumber\\
&0\leq \tau_j,\label{eq:D2p}\\
&\displaystyle\sum_{j\in J}\tau_j=1,\label{eq:D3p}\\
&s(\tau)=\displaystyle\sum_{j\in J}\tau_jx_j.\label{eq:D4p}
\end{align}

Consider an optimal solution $\tau^*$ to Problem D2. Therefore $s(\tau^*)=C$. If $\tau_j^*\geq \frac{1}{n+1}$ for all $j\in J$ then $\tau^*$ must also be optimal for Problem D1$'$. By complementary slackness the fact that $\tau_j^*\geq \frac{1}{n+1}>0$ implies that $\|C-x_j\|-\displaystyle\max_{i\in J}\|C-x_i\|=0$. Hence all points of $X$ are concyclic about $C$. To see what the domain restriction $\tau_j^*\geq \frac{1}{n+1}$ means geometrically we first need a lemma.

\begin{lem}If there exists affine weights $\tau_j\geq \frac{1}{n+1}$ such that $s=\sum\tau_jx_j$ then $s\in \mathrm{conv}(\mathcal{M})$.
\end{lem}
\begin{pf}
If $\tau_j\geq \frac{1}{n+1}$ for all $j\in J$ then $s=\sum\frac{1}{n+1}x_j+\sum u_jx_j$ where $u_j\geq 0$ and $\sum u_j=\frac{1}{n+1}$. Let $y=\dfrac{\sum u_jx_j}{\sum u_j}=(n+1)\sum u_jx_j$; note that $y\in\mathrm{conv}(X)$. Hence $s=\frac{n}{n+1}M+\frac{1}{n+1}y$, which implies that $\|s-M\|= \frac{1}{n+1}\|M-y\|\leq \frac{1}{n+1}\|M-y'\|$, where $y'$ is on the boundary of $\mathrm{conv}(X)$ and $M,y,y'$ are collinear. The result now follows from Lemma \ref{divrat}.
\end{pf}

The converse of this lemma is not true in general, however, by similar reasoning we have:
\begin{lem}If $s\in \mathrm{conv}(\mathcal{M})$ and $s$ is equidistant from all nodes of $X$ then there exists affine weights $\tau_j\geq \frac{1}{n+1}$ such that $s=\sum\tau_jx_j$.
\end{lem}

Therefore we have the following final result on $C$, which follows directly from the above discussion and lemma.

\begin{props}\label{piProp}Suppose that $C$ is equidistant from all points of $X$. Then $C$ is the min-power centre of $X$ if and only if $C\in \mathrm{conv}(\mathcal{M})$.
\end{props}

In some sense it seems we cannot say anything about the location of $s^*$ without involving $M$. In fact, let $Q$ be a polygon with at most $n$ extreme points. For any point set $X$ with convex hull $Q$, if there is no restriction on the location of $M$ then there is no restriction on the location of $s^*$ (except, of course, that $s^*$ must lie in $Q$). We make this notion clear in the next proposition before closing the section with a proof that the centroid is a good approximation for the min-power centre.

\begin{figure}[htb]
  \begin{center}
    \includegraphics[scale=0.7]{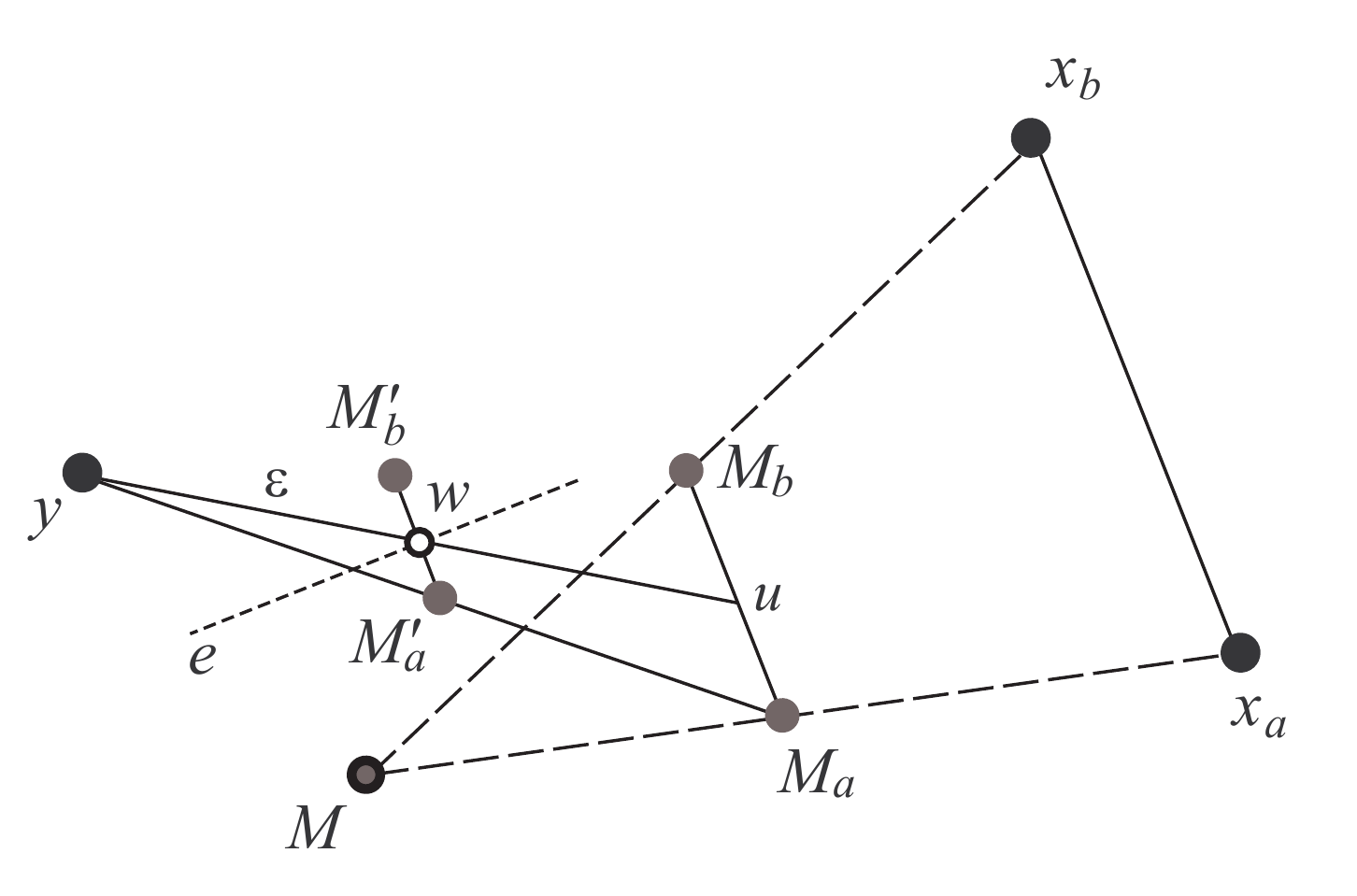}\\
  \end{center}
  \caption{Constructing an optimal point at $w$ by introducing $m$ new nodes at point $y$}
  \label{anyPoint}
\end{figure}

\begin{props}For any point $w\in \mathrm{int}(\mathrm{conv}(X))$, where $w$ is on a vertex or edge of $\mathcal{V}$, there exists a set of points $X'$ such that the extreme points of $\mathrm{conv}(X)$ and $\mathrm{conv}(X')$ coincide and such that the min-power centre of $X'$ is $w$.
\end{props}
\begin{pf}
We demonstrate the case when $w$ is in the interior of an edge $e$ of $\mathcal{V}$; the other case is similar. Let $X'$ initially consist of the extreme points of $\mathrm{conv}(X)$ only. We introduce points to the interior of $\mathrm{conv}(X')$ to get the final $X'$, which means that $\mathcal{V}(X)=\mathcal{V}(X')$ throughout. Suppose that $e$ is a portion of $V(x_a)\cap V(x_b)$ for some $a,b$. Then $w$ is the min-power centre of $X'$ if and only if $w$ also lies on the segment $M_aM_b$ -- therefore suppose that this is not the case. Let $u=\dfrac{1}{2}(M_a+M_b)$ and let $y\in\mathrm{int}(\mathrm{conv}(X))$ be a point on the line intersecting $u$ and $w$ such that $w$ is between $y$ and $u$; see Figure \ref{anyPoint}. Let $\varepsilon=\|y-w\|>0$.We introduce $m>0$ new nodes at the point $y$ and, for $j=a,b$, let $M_j'=\frac{1}{n+m+1}\left(x_j+my+\sum_{i\in J}x_i\right)$. Clearly $M_j'$ is the centre of mass of $M_j$ and $y$ (weighted $n+1$ and $m$ respectively), and therefore $M_j'$ divides the segment $yM_a$ into the ratio $\|M_a-M_a'\|:\|M_a'-y\|=m:n+1$; in fact, the segments $x_bx_a$, $M_bM_a$ and $M_b'M_a'$ are all parallel. In order for $w$ to be optimal for the resultant $X'$ the segment $M_a'M_b'$ must intersect $w$, which is equivalent to requiring that $\varepsilon =\frac{1}{m}(n+1)\|w-u\|$. The integer $m$ can therefore be selected to be large enough so that the point $y$ resulting from the subsequent value of $\varepsilon$ will be in the convex hull of $X$.
\end{pf}

In the case when $w$ from the previous proposition lies in the interior of a region of $\mathcal{V}$, a similar process will produce an optimal point arbitrarily close to $w$. It is also relatively straightforward to extend the proposition to the case where no two nodes may coincide, and to the case where all nodes of $X'$ are required to lie on the boundary of the convex hull of $X$. Any other interesting general restrictions on the location of $s^*$ that do not involve $M$ would therefore probably be for sets of nodes in convex position (i.e., where every node is an extreme point). We leave this as an open question.

The centroid is very easy to calculate -- it can be done in constant time given the locations of the $x_j$. Corollary \ref{thisCor} shows that $M$ is an appropriate point for approximating $s^*$, especially in applications with large $n$. For any set of points $X$ let $$\rho(X)=\dfrac{P(M)}{P(s^*)}.$$ To place an upper bound on $\rho$ we first need a lemma. Recall once again that $x_r$ is the farthest node from $M$.

\begin{lem}\label{mgeo}$\|M_r-x_r\|\leq \|s^*-x_r\|$.
\end{lem}
\begin{pf}
\begin{alignat*}{5}
&\|x_r-M\|& &\leq\|s^*-x_r\|+\|s^*-M\|&\quad &{(\mathrm{by\ the\ triangle\ inequality})}&\\
&{}& &\leq\|s^*-x_r\|+\|M_r-M\|&\quad &{(M_r\mathrm{\ is\ the\ farthest\ point\ of\ }}&\\
&{}&&{}&\quad &\mathrm{\ \ conv}(\mathcal{M})\mathrm{\ from\ }M).&
\end{alignat*}
Therefore
\begin{alignat*}{5}
&\|s^*-x_r\|& &\geq\|x_r-M\|-\|M_r-M\|&\quad &{}&\\
&{}& &=\|M_r-x_r\|.&\quad &{}&
\end{alignat*}
\end{pf}

\begin{thm}\label{mainRat}Let $k=|\underset{{j}}{\mathrm{arg\,max}}\,\, \{\|s^*-x_j\|\}|$. Then $$\rho(X)\leq \dfrac{1}{k+1}\left(\dfrac{n+1}{n}\right)^2+\dfrac{k}{k+1}.$$
\end{thm}
\begin{pf}
Recall that for any point $y$ we denote by $F(y)$ the sum $F(y)=\displaystyle\sum_{j\in J}\|y-x_i\|^2$. Let $x_p\in X$ be a farthest node from $s^*$. By the triangle inequality $\|M-x_r\|\leq \|s^*-M\|+\|s^*-x_r\|\leq \|s^*-M\|+\|s^*-x_p\|$.

Therefore
\begin{alignat}{4}
&\dfrac{\|M-x_r\|}{\|s^*-x_p\|}&&\leq\dfrac{\|s^*-M\|}{\|s^*-x_p\|}+1&\nonumber\qquad &{}&\\
&{}&&\leq \dfrac{\|M_r-M\|}{\|s^*-x_r\|}+1&\qquad &(\mathrm{since\ } s^*\in\mathcal{M})&\nonumber\\
&{}&&\leq \dfrac{\|M_r-M\|}{\|M_r-x_r\|}+1&\qquad &(\mathrm{by\ Lemma\ \ref{mgeo}})&\nonumber\\
&{}&&=\dfrac{1}{n}+1=\dfrac{n+1}{n}&\qquad &(\mathrm{by\ Lemma\ \ref{divrat}}).&\nonumber
\end{alignat}

We now have
\begin{alignat}{4}
&\dfrac{P(M)}{P(s^*)}& &\leq \dfrac{\|M-x_r\|^2+F(s^*)}{\|s^*-x_p\|^2+F(s^*)}& &(\mathrm{since\ }M\mathrm{\ minimises\ }F)&\nonumber\\
&{}& &=\dfrac{\|M-x_r\|^2-\|s^*-x_p\|^2}{\|s^*-x_p\|^2+F(s^*)}+1& &{}&\nonumber\\
&{}& &\leq\dfrac{\|M-x_r\|^2-\|s^*-x_p\|^2}{(k+1)\|s^*-x_p\|^2}+1&\qquad &(\mathrm{since\ }k\mathrm{\ longest\ edges\ are\ equal})&\nonumber\\
&{}& &=\dfrac{\|M-x_r\|^2}{(k+1)\|s^*-x_p\|^2}+\dfrac{k}{k+1}& &{}&\nonumber\\
&{}& &\leq \dfrac{1}{k+1}\left(\dfrac{n+1}{n}\right)^2+\dfrac{k}{k+1}.& &{}&\nonumber
\end{alignat}
\end{pf}

Since $k\geq 1$ we conclude that $\rho(X)\leq \dfrac{1}{2}\left(\dfrac{n+1}{n}\right)^2+\dfrac{1}{2}$ for any $X$. The above approximation is ``good" in the following sense. Let $X_n$ be any set of $n$ points.

\begin{cors}\label{thisCor}$\displaystyle\lim_{n \rightarrow \infty}\rho(X_n)=1.$
\end{cors}

\section{A brief look at general $\alpha>1$}\label{genA}
The situation when $\alpha\neq 2$ is complicated by the fact that the level curves of the $F_j$ are no longer circular. In particular, this implies that there are no easily constructed and useful fixed points directly analogous to the $M_j$. Even though the farthest point Voronoi diagram still plays a crucial role for general $\alpha$, there is no simple dual structure analogous to the Delaunay triangulation. We may consider $\underset{{s}}{\mathrm{arg\,min}}\,\, \sum\|s-x_j\|^\alpha$ (the generalised Fermat-Weber point) as an analogue for $M$, however, a compass and straight-edge construction of this point is almost certainly impossible since even the classical Fermat-Weber point is not constructible in this way (see \cite{baj}).

This section will show, however, that not all insights gained for the $\alpha=2$ case are lost when $\alpha\neq 2$. In fact, there are several important uses for the mathematical machinery developed in the preceding sections. Firstly, as discussed in the introduction, path loss exponents close to $\alpha=2$ often occur in real radio network scenarios. It seems evident that $s^*$ (and even $M$) will serve as a good approximation points for min-power centres based on $\alpha$ close to $2$.

Also, some of the the results from the previous section, especially those involving the $1$-centre $C$, are directly applicable for any $\alpha>1$. The conceptual framework that we present for proving these results consists of transforming the set $X$ into a new set $X(s)$ such that the general $\alpha>1$ min-power centre problem on $X$ is transformed into the quadratic min-power centre problem on $X(s)$. Preliminary research indicates that this framework will also be useful for the development of bounds on the performance of $s^*$ or the centroid as approximation points for general $\alpha$.

\begin{figure}[htb]
  \begin{center}
    \includegraphics[scale=0.45]{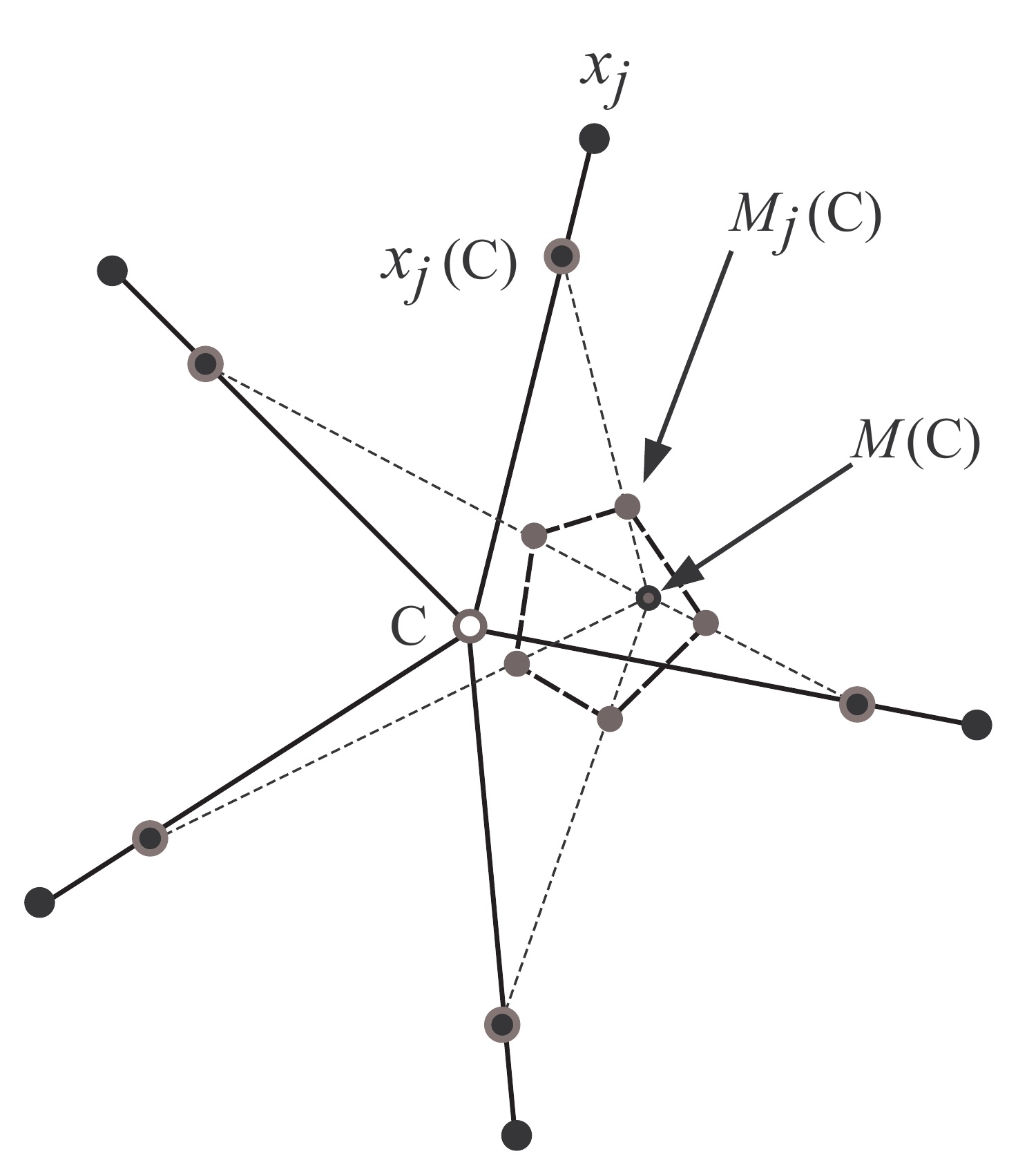}\\
  \end{center}
  \caption{An example for Proposition \ref{toComb} when $\alpha<2$ and $C$ is not optimal}
  \label{figPiZero}
\end{figure}

We first state some notation and a lemma. Recall, from Section \ref{SecProps}, that the power of the star $G$ centred at $s$ with leaf-set $X$ and for a given $\alpha$ is denoted $P_\alpha(s,X)$, and that $P_\alpha(s,X)=\displaystyle\max_{j\in J}\{F_j^\alpha(s,X)\}$. We also similarly generalise the definition of $F(s)$ to $F^\alpha(s,X)$. The point minimising $P_\alpha(s,X)$ is now denoted by $s_\alpha^*$. Next we describe the construction of the set $X(s)$. Essentially, each vector $x_i-s$ is transformed into a vector $x_i(s)-s$ such that $\|x_i-s\|^{\alpha-1}=\|x_i(s)-s\|$. The convex hull of the resultant set of nodes will therefore contain the convex hull of $X$ when $\alpha\geq 2$, with the inverse of this property holding when $\alpha<2$. In formal terms, for any point $s$ let $x_i(s)=s+\|s-x_i\|^{\alpha-2}(x_i-s)$ and let $X(s)=\{x_j(s):j\in J\}$. Let $M(s)$ be the centre of mass of $X(s)$, and let $M_j(s)$ be the $j$-th $2$-centroid of $X(s)$. Some of these concepts are illustrated in Figure \ref{figPiZero} for the set $X(s)$ with $s=C$ and $\alpha<2$.

\begin{lem}\label{lemS}Let $s'$ be any point. Then $s=s'$ minimises $P_2(s,X(s'))$ if and only if $s'=s_\alpha^*$.
\end{lem}
\begin{pf}
Suppose that $s'=s_\alpha^*$. Since $s_\alpha^*$ is optimal for $P_\alpha$ the following KKT conditions are satisfied: there exist multipliers $\tau_j^*$ such that (A): $\dfrac{1}{n+1}\leq \tau_j^*$ for all $j$; (B): $\sum\tau_j^*=1$; (C): $\sum\tau_j^*(s_\alpha^*-x_j)\|s_\alpha^*-x_j\|^{\alpha-2}=0$; and (D): $\tau_j^*\left(\|s_\alpha^*-x_j\|-\displaystyle\max_i\{\|s_\alpha^*-x_i\|\}\right)=0$ for all $j$. Note that (C) can be written as $\sum\tau_j^*(s_\alpha^*-x_j(s_\alpha^*))=0$. Therefore conditions (A)-(C) are the first three KKT conditions for $s_\alpha^*$ being optimal for $P_2$ on $X(s_\alpha^*)$. We only need to still show complementary slackness. Note that $\|s_\alpha^*-x_j(s_\alpha^*)\|=\|s_\alpha^*-x_j\|^{\alpha-1}$, so that if $\tau_j^*>0$ then $\|s_\alpha^*-x_j\|=\displaystyle\max_i\{\|s_\alpha^*-x_i\|\}$, which implies that $\|s_\alpha^*-x_j(s_\alpha^*)\|=\displaystyle\max_i\{\|s_\alpha^*-x_i(s_\alpha^*)\|\}$. Therefore complementary slackness also holds for $P_2$ at $s_\alpha^*$. The other direction follows similarly.
\end{pf}

The next proposition is a direct generalisation of Proposition \ref{piProp} and follows from the previous lemma. An example is illustrated in Figure \ref{figPiZero}.

\begin{props}\label{toComb}Suppose that $C$ is equidistant from all points of $X$. Then $s_\alpha^*=C$ if and only if $C\in \mathrm{conv}(\{M_j(C)\})$.
\end{props}
\begin{pf}
Suppose that $s_\alpha^*=C$. Then, since $s_\alpha^*$ minimises $P_2(s,X(s_\alpha^*))$, we must have $C=s_\alpha^*\in \mathrm{conv}(\{M_j(s_\alpha^*)\})=\mathrm{conv}(\{M_j(C)\})$ by KKT Condition (\ref{eq:KKT1}). Next suppose that $C\in \mathrm{conv}(\{M_j(C)\})$. Note that $C$ is also the $1$-centre of $X(C)$ since all points are equidistant from $C$. Therefore $C$ minimises $P_2(s,X(C))$ by Proposition \ref{piProp}, and therefore minimises $P_\alpha(s,X)$ by Lemma \ref{lemS}.
\end{pf}

This proposition can be strengthened, as shown the following result.

\begin{props}\label{equiS}If $s_2^*$ is equidistant from all point of $X$ then $s_\alpha^*=s_2^*$.
\end{props}
\begin{pf}
Observe that $X(s_2^*)$ is a uniform scaling of $X$ from the point $s_2^*$. Since $s_2^*$ minimises $P_2(s,X)$ it also minimises $P_2(s,X(s_2^*))$. Therefore, by Lemma \ref{lemS}, $s_2^*=s_\alpha^*$.
\end{pf}

As a final observation note that as $\alpha$ tends to infinity the longest edge incident to $s_\alpha^*$ dominates $F^\alpha(s,X)$. Therefore, since $C$ minimises $\max\{\|s-x_j\|^\alpha\}$, and $P_\alpha(s,X)=\max\{\|s-x_j\|^\alpha\}+F^\alpha(s,X)$, we have the following result.

\begin{observation}$\displaystyle\lim_{\alpha \rightarrow \infty}s_\alpha^*=C.$
\end{observation}

\section{Conclusion}
In this paper we solve the quadratic min-power centre problem in the Euclidean plane. We provide a complete geometric description and method of construction for the optimal point by means of farthest point Voronoi diagrams. The solution leads to various structural results relating the min-power centre to the centroid and the $1$-centre, and, in particular, allows us to construct an explicit bound on the performance of the centroid as an approximation point. We anticipate that the mathematical machinery developed in this paper will be applied to an important fundamental model for the physical design of wireless ad hoc networks -- namely the geometric range assignment problem where a bounded number of new nodes may be introduced anywhere in the plane. Preliminary research has demonstrated that the framework we develop in this paper will also be useful for problems with cost functions that are not quadratic.

\end{document}